\documentclass[12pt,centertags,oneside]{amsart}
\usepackage{amsmath,amstext,amsthm,amscd,typearea}
\usepackage{amssymb}
\usepackage{a4wide}
\usepackage[mathscr]{eucal}
\usepackage{mathrsfs}
\usepackage{typearea}
\usepackage{charter}
\usepackage{pdfsync}
\usepackage{url}

\usepackage[a4paper,width=16.2cm,top=3cm,bottom=3cm]{geometry}

\numberwithin{equation}{section}

\newtheorem{theorem}{Theorem}[section]

\newtheorem{proposition}[theorem]{Proposition}

\newtheorem{lemma}[theorem]{Lemma}

\newcommand{\cali}[1]{\mathscr{#1}}

\newcommand{\vol}{\mathop{\mathrm{vol}}}

\newcommand{\ddc}{\text{\normalfont dd}^c}
\newcommand{\dc}{\text{\normalfont d}^c}

\newcommand{\capK}{\rm cap}

\newcommand{\B}{\mathbb{B}}
\newcommand{\C}{\mathbb{C}}

\newcommand{\N}{\mathbb{N}}

\title[]{Locally pluripolar sets are pluripolar}

\author{Duc-Viet Vu}
\address{University of Cologne, Mathematical Institute, Germany}
\email{vuduc@math.uni-koeln.de}

\date{\today}
\begin{document}

\begin{abstract}  We prove that every locally pluripolar set on a compact complex manifold is pluripolar.  This extends similar results in K\"ahler case. 
\end{abstract}

\maketitle

\medskip

\noindent
{\bf Classification AMS 2010}:  32W20, 32U40.  

\medskip



\section{Introduction} \label{introduction}

Pluripotential theory has been a crucial tool in complex geometry, complex dynamics as well as other fields of Mathematics. We refer to \cite{Demailly_ag,DS_book,Kolodziej05,Levenberg} for some expositions of this theory and its applications. Among other things, locally pluripolar sets are important objects in the pluripotential theory which play the role of negligible sets as a counterpart to  algebraic subvarieties in algebraic geometry, see the next section for definitions. To illustrate this comparison, we recall that locally pluripolar sets are of Hausdorff codimension at least $2$ (see \cite[Th. 3.13]{Landkof}) and their intersections with totally real submanifolds of the ambient manifold are of Lebesgue measure zero (see \cite[Cor. 1.2]{Vu_MA}).  We refer to \cite{Klimek,Skoda_integrability,DVS_exponential,Lucas_skodatype} for more information.   

Josefson's theorem \cite{Josefson}, which is a key result in the pluripotential theory on $\C^k,$ affirms that locally pluripolar sets on $\C^k$ are in fact (globally) pluripolar. Simplified proofs of this fact  were given by Bedford-Taylor \cite{Bedford_Taylor_82} and  Alexander-Taylor \cite{ATaylor_capacity}.   This result was generalized to the pluripolar sets on projective manifolds, compact K\"ahler manifolds, respectively,  by Dinh-Sibony \cite{DS_tm}, Guedj-Zeriahi  \cite{GZ}, see also Berman-Boucksom-Witt Nystr\"om \cite{BoucksomBermanWitt} for the case of manifolds equipped a big line bundle.   Our main result below extends this property to pluripolar sets on \emph{every} compact complex manifolds. 
 
\begin{theorem} \label{the_pluripolar} Every locally pluripolar set on a compact complex manifold is pluripolar.
\end{theorem}

By the above theorem,   there exist abundantly  non-continuous quasi-p.s.h. functions on $X.$ This is a fact which probably cannot be seen directly because unlike projective manifolds, a general compact complex manifold might have very few hypersurfaces. The key ingredients of the proof of Theorem \ref{the_pluripolar} are the comparison  (\ref{ine_TAorigin}) between capacities  generalizing similar comparison results in \cite{ATaylor_capacity, GZ} and recent developments of the pluripotential theory for non-K\"ahler manifolds by Ko{\l}odziej, Dinew and  Nguyen \cite{KD_hermitian,KC_hermitian}.   \\

\noindent
\textbf{Acknowledgments.} The author would like to thank Ngoc Cuong Nguyen for fruitful comments.  This research is supported by a postdoctoral fellowship of  the Alexander von Humboldt Foundation. 

\section{Proof of Theorem \ref{the_pluripolar}}

First of all, we need to recall some basic notations from the pluripotential theory. Let $X$ be a compact complex manifold of dimension $k$. A function from $X$ to $[-\infty, \infty)$ is said to be  \emph{quasi-p.s.h.} if it can be written  locally as  the sum of a plurisubharmonic (p.s.h.) function and a smooth one.   Put $\dc:= i/(2\pi)(\overline \partial -\partial).$  For a continuous real  $(1,1)$-form $\eta$,  a quasi-p.s.h. function $\varphi$ is said to be \emph{$\eta$-p.s.h.} if $\ddc \varphi+ \eta\ge 0.$ We have the following characterization of quasi-p.s.h. functions in terms of  submean-type inequalities. 

\begin{lemma} \label{le_charac} Let $U$ be an open subset of $\C^k$ and $\eta$ a continuous real $(1,1)$-form  on $U.$ A  function $\varphi: U \to [-\infty, \infty)$  is  $\eta$-p.s.h. if and only if it is upper semi-continuous, not identically $-\infty$ and  for every $x \in U$ and  every complex line $L_v:=\{x+ t v : t \in \C\},$ for some $v \in  \C^k,$ passing through $x,$ we have
\begin{align} \label{ine_submean}
\varphi(x) \le  \frac{1}{2\pi} \int_{0}^{2 \pi} \varphi(x+\epsilon  e^{i \theta}v) d\theta + \int_0^\epsilon \frac{d t}{t} \int_{\{|s|\le t\}} \eta_v,
\end{align}
for every constant $\epsilon>0$ small enough, where $\eta_v(t)$ is the restriction of $\eta$ to  $L_v$ which is identified with $\C$ via $t \longmapsto x+ tv.$ 
\end{lemma}

\proof Consider an $\eta$-p.s.h. function $\varphi.$ We need to verify (\ref{ine_submean}).   For every positive constant $r,$ let $\chi_r$ be a smooth multi-radial nonnegative function  compactly supported on the polydisk of radius $r$ in $\C^k$ with $\int_{\C^k} \chi_r(x) \vol(x)=1,$ where $\vol$ is the canonical volume form on $\C^k.$ Since $\varphi$ is locally integrable, we can define the convolution
 $$\varphi^r(x):= \int_{\C^k} \varphi(x-y) \chi_r(y) \vol(y)$$
which is smooth. We have $\varphi^r \to \varphi$ pointwise  as $r\to 0$ because $\varphi$ can be written as the sum of a p.s.h. function and a smooth one.  Denote by 
$$\eta^r(x):=\int_{\C^k} \eta(x-y)\chi_r(y) \vol(y)$$
which converges uniformly to $\eta$  as $r \to 0$ because $\eta$ is continuous.   Hence,  $\ddc \varphi +\eta \ge 0$ if and only if $\ddc \varphi^r+ \eta^r \ge 0$ for every $r$ small. Similarly, (\ref{ine_submean}) holds if  it holds for $(\varphi^r, \eta^r)$ in place of $(\varphi, \eta)$ for every small $r.$ It follows that it suffices to prove (\ref{ine_submean}) for smooth $\varphi$ and smooth $\eta.$ 

Hence we can assume $\varphi, \eta$ are smooth and  follow standard arguments in \cite{Hormander}.  Let $v \in \C^k$ and $x \in U.$  
Put $\varphi_v(t):= \varphi(x+ t v).$  We get $\ddc \varphi_v + \eta_v \ge 0.$ The Lelong-Jensen formula for $\varphi_v(t)$ gives
$$M_{\epsilon,v} - M_{\epsilon',v}= \int_{\epsilon'}^\epsilon \frac{dt}{t} \int_{\{|s| \le t\}} \ddc \varphi_v,$$
where  $\epsilon> \epsilon'$ are positive constants and 
$$M_{s,v}:=  \frac{1}{2\pi} \int_{0}^{2 \pi} \varphi_v(\epsilon  e^{i \theta}) d\theta$$
for every constant $s>0.$ It follows that 
$$M_{\epsilon',v} \le M_{\epsilon,v} +  \int_{\epsilon'}^\epsilon \frac{dt}{t} \int_{\{|s| \le t\}} \eta_v.$$
Letting $\epsilon'\to 0$ in the last inequality
 gives (\ref{ine_submean}) because $\varphi_v$ is continuous at $0$. 
 
Assume now (\ref{ine_submean}). This combined with the hypothesis that $\varphi \not = -\infty$ implies $\varphi \in L^1_{loc}.$ Moreover, as in the case of p.s.h. functions, since $\varphi$ is upper semi-continuous,   (\ref{ine_submean}) also tells us that $\varphi$ is strongly semi-continuous in the sense that for every Borel subset $A$ of $U$ whose complement in $U$ is of zero Lebesgue measure, we have 
\begin{align} \label{eq_stronglypsh}
\limsup_{y \in A\to x} \varphi(y)= \varphi(x).
\end{align}
Consider first the case where $\varphi \in \cali{C}^2.$  Direct computations show
$$ \epsilon^{-2}\big(M_{\epsilon, v} - \varphi_v(0)\big) \to \pi \ddc \varphi_v(0)/2$$
as $\epsilon \to 0.$ Applying this to (\ref{ine_submean}) gives $\ddc \varphi_v(0) + \eta_v(0) \ge 0.$  In other words, we get $\ddc \varphi + \eta \ge 0. $   

In general, let $\varphi^r, \eta^r$ be as above. Since $\varphi \in  L^1_{loc},$ $\varphi^r \to \varphi$ in $L^1_{loc}.$  We see easily that (\ref{ine_submean}) also holds for $(\varphi^r, \eta^r)$ in place of $(\varphi, \eta).$ By the above arguments, $\ddc \varphi^r + \eta^r \ge 0.$ Letting $r\to 0$ gives $\ddc \varphi+ \eta \ge 0.$ 

It remains to check that $\varphi$ is the sum of a p.s.h. function and a smooth one. To this end, we only need to work locally. Thus, we can assume there is a smooth function  $\psi$ on $U$ with $\ddc \psi \ge \eta.$ We deduce $\ddc \varphi_1 \ge 0$ for $\varphi_1:= \varphi+ \psi$ which is also strongly semi-upper continuous in the above sense. Let $\varphi_1^r$ be the regularisation of $\varphi_1$ defined in the same way as $\varphi^r.$ Notice that $\varphi_1^r \to \varphi_1$ in $L^1_{loc}$ and $\varphi_1^r$ is p.s.h. and decreasing to some p.s.h. function $\varphi'_1.$ Hence, $\varphi_1= \varphi_1'$ almost everywhere. Using this and (\ref{eq_stronglypsh}) yield that $\varphi_1= \varphi'_1$ everywhere. In other words, $\varphi$ is quasi-p.s.h..   This ends the proof.
\endproof

The following extension result generalizes the similar property for p.s.h. functions.  

\begin{lemma} \label{le_extensionpsh} Let $U$ be an open subset in a complex manifold $Y$. Let $\eta$ be a continuous  real $(1,1)$-form  on $Y.$  Let $\psi_1$ be an $\eta$-p.s.h. function on $U$ and  $\psi_2$  an $\eta$-p.s.h function on $Y$ such that $\limsup_{y \to x}\psi_1(y) \le  \psi_2(x)$ for every $x \in \partial U$. Define $\psi:= \max\{\psi_1, \psi_2\}$ on $U$ and $\psi:= \psi_2$ on $Y \backslash U.$ Then $\psi$ is an $\eta$-p.s.h. function.
\end{lemma}

\proof This is a direct consequence of Lemma \ref{le_charac}.
\endproof

  A subset $A$ of  $X$ is  \emph{locally pluripolar} if  every point $x$ in $A$ there is an open neighborhood $U_x$ of $x$ in $X$  and a p.s.h. function $\varphi$ on $U_x$ for which $A \cap U_x \subset \{\varphi=-\infty\}.$   A subset $A$ of $X$ is \emph{pluripolar} if $A \subset \{\varphi= -\infty\}$ for some quasi-p.s.h. function $\varphi$ in $X.$   

For every Borel set $A'$ in an open subset $U$ of $\C^k$,  Bedford-Taylor \cite{Bedford_Taylor_82} introduced the following notion of   \emph{capacity} of $A'$ in $U:$  
$$\capK_{BT}(A', U):= \sup \big\{ \int_{A} (\ddc \varphi)^k: \varphi \, \text{ p.s.h.  on $U$ }, 0 \le \varphi \le 1 \, \text{on U}\big\}.$$
Fix, from now on,  a Hermitian metric $\omega$ on  $X.$    For every Borel set $A\subset X,$ define
$$\capK_{BTK}(A):= \sup \big\{ \int_{A} (\ddc \varphi+\omega)^k: \varphi \, \text{ $\omega $-p.s.h.}, 0 \le \varphi \le 1\, \text{on X}\big\}.$$
The last capacity was introduced by Ko{\l}odziej \cite{Kolodziej05} as an analogue to the local capacity $\capK_{BT}$    and is used  to study complex Monge-Amp\`ere equations on Hermitian manifolds, see  for example \cite{Kolodziej05,KD_hermitian,KC_hermitian,NgocCuong_hermitian}. By Lemma \ref{le_boundecap} below, $\capK_{BTK}(A)$ is always finite. It is also clear that if we use another Hermitian metric to define $\capK_{BTK}$, then the resulted capacity is equivalent to that associated to $\omega.$ 

We will need the following modified version of the classical Bedford-Taylor comparison principle due to Ko{\l}odziej and Nguyen, see \cite{KD_hermitian} for a related result.  

\begin{proposition}\label{pro_comparison} \cite[Th. 0.2]{KC_hermitian} Let $\varphi, \psi$ be  bounded $\omega$-p.s.h  functions on $X.$ Let $0<\epsilon<1$ and $m_\epsilon:= \inf_X(\varphi - (1-\epsilon)\psi).$ Then there exists a big constant $B>0$ depending only on $\omega, k$ such that   for every constant $0< s < \epsilon^3/(16 B)$ we have
$$\int_{\{\varphi < (1-\epsilon) \psi+m_\epsilon+s \}} \big((1-\epsilon)\ddc \psi+ \omega\big)^k \le (1+ C \epsilon^{-k}s) \int_{\{\varphi < (1-\epsilon) \psi+m_\epsilon+s \}} (\ddc \varphi+ \omega)^k,$$
where $C$ is a constant depending only on $k,B.$ 
\end{proposition}

A consequence of the last result is the following.  

\begin{lemma} \label{le_boundecap} (\cite{KD_hermitian, KC_hermitian}) Let $M$ be a positive number. Then there exists a constant $c_M>0$ such that for every $\omega$-p.s.h. function $\varphi$ bounded by $M,$ we have  
\begin{align} \label{ine_boundeMAhermia}
0<  \int_X (\ddc \varphi+ \omega)^k \le c_M.
\end{align} 
\end{lemma}
\noindent
However, we don't know  whether 
$$\inf_{\{\varphi:\,  |\varphi| \le M\}} \int_X (\ddc \varphi+ \omega)^k >0?$$
\proof   
 The  second desired inequality is proved in \cite{KD_hermitian} by using integration by parts. The first one is  observed in \cite{KC_hermitian}. To see it, it is enough to notice that by choosing $\epsilon:= 1/2$ and $s>0$ small enough in Proposition \ref{pro_comparison},  for every $\omega$-p.s.h. $\psi$ with $0 \le \psi \le s$ and $\varphi$ as in the hypothesis, we have
$$\int_{\{\varphi < \inf_X \varphi+ s\}} \big(\ddc \psi+ \omega\big)^k \lesssim \int_{\{\varphi < (1-\epsilon) \psi+m_\epsilon+2s \}} \big((1-\epsilon)\ddc \psi+ \omega\big)^k \lesssim_s \int_{X} (\ddc \varphi+ \omega)^k$$
because 
$$\{\varphi < \inf_X \varphi+ s\}  \subset \{\varphi < (1-\epsilon) \psi+m_\epsilon+ 2s \}.$$
It follows that  
\begin{align} \label{cor_modificat}
\int_X (\ddc \varphi+ \omega)^k \gtrsim_s \capK_{BTK}\big(\{\varphi < \inf_X \varphi+s\}\big)
\end{align}
which is strictly positive because it is the capacity of a non-empty open set. The proof is finished. 
\endproof

Let $(U_j)_{1 \le j \le N}$ and $(U'_j)_{1 \le j \le N}$ be  finite open coverings of $X$ such that  $\overline U_j$ is smooth and  contained in some local chart of $X$ biholomorphic to a polydiscs  for every $1 \le j \le N,$   $U_j= \{\psi_j <0 \}$ for some p.s.h. function $\psi_j$ defined on an open neighborhood of $\overline U_j$  with $\partial U_j = \{\psi_j=0\}$ and $U'_j \Subset U_j$ for $1 \le j \le N.$ In practice, it suffices to take $U_j, U_j'$ to be balls and $\psi_j$ are the differences of radius functions and  constants. 

\begin{lemma} \label{le_sosanh2capa} (\cite{Kolodziej05,KD_hermitian}) There exists strictly positive constants  $c_1, c_2$ such that  for every $A \subset X$ we have
$$c_1 \sum_{j=1}^N \capK_{BT}\big(A\cap U'_j, U_j\big) \le  \capK_{BTK}(A) \le c_2 \sum_{j=1}^N \capK_{BT}\big(A\cap U'_j, U_j\big).$$
\end{lemma} 

\proof Put $A'_j:= A \cap U'_j$ which is a relatively compact subset of $U_j.$ We have $\cup_j A'_j= A.$ The second desired inequality is obvious from the definitions of capacities. We prove now the first desired inequality.   

Fix an index $1 \le j \le N.$ By our choice of $U_j,$ for  every p.s.h. function $0 \le u \le 1$ on $U_j,$ we can find another p.s.h. function $-1 \le \tilde{u} \le 0$ on $U_j$ satisfying $\tilde{u}= u-1$ on some open neighborhood of $\overline U'_j$ and $\tilde{u}= 0$ on $\partial U_j.$ Such a  $\tilde{u}$ can be chosen to be $\max\{u-1,  A \psi_j\}$ for  some constant $A$ big enough.  Clearly, 
$$\int_{A'_j} (\ddc u)^k= \int_{A'_j} (\ddc \tilde{u})^k.$$    
Since $-1 \le \tilde{u} \le 0$ and $\tilde{u}= 0$ on $\partial U_j,$  there is a quasi-p.s.h. function $\tilde{u}_1$ on $X$ such that $\ddc \tilde{u}_1 + C \omega \ge 0$ for some constant $C$ independent of $\tilde{u}$ and $\tilde{u}_1=\tilde{u}$ on some open neighborhood of $\overline U'_j.$  We deduce that 
$$\int_{A'_j} (\ddc u)^k = \int_{A'_j} (\ddc \tilde{u}_1)^k \le \int_{A'_j} (\ddc \tilde{u}_1+ C \omega)^k \le C^k \capK_{BTK}(A'_j).$$
Consequently, $\capK_{BT}(A'_j, U_j) \le C^k \capK_{BTK}(A'_j)$. Summing over $1\le j \le N$ in the last inequality gives the first desired inequality. This finishes the proof.       
\endproof

Since we already know that if $A$ is locally pluripolar in $U$, then $\capK_{BT}(A,U)=0$ (see  \cite[Th. 4.6.4]{Klimek} or \cite{Bedford_Taylor_82}), we get $\capK_{BTK}(A)=0$ if $A$ is locally pluripolar in $X$.  Let $(u_j)$ be a family of p.s.h. functions on an open subset $U$ of $\C^k$ locally bounded from above. Define $u:= \sup_j  u_j$ and $u^*:= \sup^*_j u_j$ the upper semi-continuous regularisation of $u.$ The set $\{u< u^*\}$ is called \emph{a negligible set} in $U$.  By Bedford-Taylor \cite{Bedford_Taylor_82},  the negligible sets are locally pluripolar. 
The following notion of capacity, which is related to  those of Alexander \cite{Alexander} and Sibony-Wong \cite{Sibony_Wong},    is due to Dinh-Sibony \cite{DS_tm}:  for $A \subset X,$
$$\capK_{ADS}(A):= \inf\{ \exp(\sup_A \varphi): \varphi \, \text{ $\omega $-p.s.h.  on $X$},\, \sup_X \varphi=0\},$$
see  \cite{Harvey_Lawson} for some of its applications. 
\begin{lemma} \label{le_capADS} $\capK_{ADS}(A)=0$ if and only if $A$ is pluripolar on $X$. 
\end{lemma}

\proof  If $A \subset  \{\varphi =-\infty\}$ for some quasi-p.s.h. $\varphi,$ it is clear that $\capK_{ADS}(A)=0.$ Consider now 
\begin{align}\label{eq_capADSA}
\capK_{ADS}(A)=0.
\end{align}
Recall that there exists a constant $c$ such that for every $\omega$-p.s.h. function $\varphi$ with the normalization condition $\sup_X \varphi=0,$ we have 
\begin{align}\label{eq_capADSA2}
\|\varphi\|_{L^1(X)} \le c.
\end{align}
We refer to \cite{Hormander,DS_tm,KD_hermitian} for a proof.  Using (\ref{eq_capADSA}), there exists a sequence of $\omega$-p.s.h. functions $(\varphi_n)$ with $\sup_X \varphi_n =0$ such that $\sup_A \varphi_n \le -n^3.$ Put 
$$\varphi:= \sum_{n=1}^\infty \frac{\varphi_n}{n^2}$$
which is a well-defined quasi-p.s.h. function because of (\ref{eq_capADSA2}). On the other hand, 
$$\sup_A \varphi \le \sum_{n=1}^\infty \frac{-n^3}{n^2}= -\infty.$$
It means that $A\subset  \{\varphi = -\infty\}.$ This finishes the proof. 
\endproof

Let $(\varphi_j)_{j\in J}$ be a family of $\omega$-p.s.h. functions uniformly bounded from above. Define $$\varphi_J:= \sup_{j \in J} \varphi_j.$$ 
Observe that $\varphi_J^*$ is an $\omega$-p.s.h. function. This can be seen by using Lemma \ref{le_charac} or  noticing that for every $\omega$-p.s.h. functions $\varphi_j, \varphi_{j'}$ we have $\max\{\varphi_j, \varphi_{j'}\}= \lim_{n \to \infty} n^{-1} \log (e^{n \varphi_j}+ e^{n \varphi_{j'}})$ whose $\ddc$ is $\ge -\omega$ for every $n$.  As in the  local setting,  $\{\varphi^*_J> \varphi_J\}$ is a locally pluripolar set.  We will present below an important case of $(\varphi_j)_{j \in J}$ and its associated extremal function $\varphi^*_J.$ 

Let $A$ be a \emph{non-pluripolar} subset of $X.$  As in the local setting or in  the K\"ahler case, we introduce the following extremal $\omega$-p.s.h. function:
$$T_A:= \sup \big\{ \varphi  \, \text{ $\omega $-p.s.h.}:  \varphi \le 0   \text{ on } A \big\}.$$
It is clear that $T_A \ge 0.$   Let $T_A^*$ be the upper semi-continuous regularisation of $T_A$.   We can check that  
\begin{align} \label{eq_ADSTA}
\capK_{ADS}(A)= \exp(-\sup_X T_A ).
\end{align}
Thus $T_A$ is bounded from above  because  $A$ is non-pluripolar. We deduce that $T^*_A$ is a bounded $\omega$-p.s.h. function and $Q_A:=\{T_A^*> T_A\}$ is a locally pluripolar set. This combined with the fact that $T_A= 0$ on $A$ implies that  $T^*_A=0$ on $A \backslash Q_A.$  
The following generalized a well-known property of $T^*_A$ in the K\"ahler case.

\begin{proposition} \label{pro_MATAuA}  Let $A$ be a nonpluripolar compact subset of $X.$  
 We have
\begin{align}\label{eq_MATAuA}
(\ddc T^*_A+ \omega)^k= 0
\end{align}
 on $X \backslash A.$
\end{proposition}

\proof We follow the usual strategy. The  key points are the existence of solutions of the Dirichlet problems proved in \cite{KC_hermitian,GuanLi,Cherrier,CKNP} and Lemma \ref{le_extensionpsh} above. 

By Choquet's lemma, there exists an increasing sequence of $\omega$-p.s.h. function $\varphi_n$ for which $T_A^*= (\lim_{n \to \infty} \varphi_n)^*.$  
For every $\omega$-p.s.h function $\varphi$ and every positive constant $\epsilon,$ using a regularisation  of $\varphi$ (see \cite{Blocki_kolo_regular}),  Hartog's lemma and the compactness of $A$, we deduce that there exists a smooth $\omega$-p.s.h. function $\varphi'$ such that $\varphi \le \varphi'$ and $\varphi' \le  \sup_K\varphi+ \epsilon$ on $K.$  We  construct  a sequence $(\varphi''_n)$ of smooth $\omega$-p.s.h. functions  from $(\varphi_n)$ inductively as follows. Let $\varphi'_1$ be a smooth $\omega$-p.s.h. function such that $\varphi_1 \le \varphi'_1$ and $\varphi'_1 \le 1$ on $A.$ For $n \ge 2, $  let $\varphi'_n$ be a smooth $\omega$-p.s.h. function such that 
$$\max\{\varphi_n, \varphi'_{n-1} - (n-1)^{-2}\} \le \varphi'_n$$
and $\varphi'_n \le 1/n^2$ on $A.$ Put 
$$\varphi''_n:= \varphi'_n - \sum_{j=n}^\infty j^{-2}.$$
By our construction, $(\varphi''_n)$ is increasing and $\varphi''_n \le 0$ on $A$ and $\varphi''_n \ge \varphi_n - (n-1)^{-1}$ for $n \ge 2.$ We infer that 
$$T_A^*= (\lim_{n \to \infty} \varphi''_n)^*.$$
Let $\B$ be an open ball in $X \backslash A.$ 
By  \cite[Th. 4.2]{KC_hermitian}, there exists $\omega$-p.s.h. functions $u_n$ on $\B$ which is in $\cali{C}^0(\overline \B)$ for which $(\ddc u_n+ \omega)^n=0$ on $\B$ and $u_n= \varphi''_n$ on $\partial \B.$  Define  $\tilde{\varphi}''_n:= u_n$ on $\overline \B$ and $\tilde{\varphi}''_n:= \varphi''_n$ outside $\B.$  By the domination principle \cite[Cor. 3.4]{KC_hermitian}, we get $u_n \ge \varphi''_n$ on $\B$ and $u_{n+1} \ge u_n$ because $\varphi''_{n+1} \ge \varphi''_n.$   By Lemma \ref{le_extensionpsh}, $\tilde{\varphi}''_n$ is an $\omega$-p.s.h. function. 
We have obtained a sequence $(\tilde{\varphi}''_n)$ of  continuous $\omega$-p.s.h functions increasing  almost everywhere to $T^*_A.$ Hence, $$(\ddc \tilde{\varphi}''_n+ \omega)^k \to (\ddc T^*_A + \omega)^k$$ as $n \to \infty.$ We thus get $(\ddc T^*_A + \omega)^k=0$ on $\B$ for every $\B$ in $X \backslash A.$ The desired equality follows. This finishes the proof.
\endproof

\begin{proposition} \label{pro_sosanh2cap} Let $A$ be a nonpluripolar  compact subset of $X.$ Then there exist strictly positive constants $c_1, c_2, \lambda_1, \lambda_2$ independent of $A$ such that 
\begin{align} \label{ine_TAorigin}
\exp\big(- \lambda_1 \capK_{BTK}^{-1}(A)\big) \le \capK_{ADS}(A) \le c_2  \exp\big(- \lambda_2  M_A^{1/k}  \capK_{BTK}^{-1/k}(A)\big).
\end{align}
where $M_A:= \int_X (\ddc T^*_A+ \omega)^k>0.$ 
\end{proposition}
Note that $M_A>0$ because of Lemma \ref{le_boundecap}. 
\proof  
Since $A$ non-pluripolar, $T^*_A$ is a bounded $\omega$-p.s.h. function. By (\ref{eq_ADSTA}), the desired inequalities are equivalent to the following:
\begin{align}\label{ine_TA}
\lambda_1 \capK_{BTK}^{-1}(A)  \ge  \sup_X T_A \ge c'_2+  \lambda_2 M_A^{1/k}\capK_{BTK}^{-1/k}(A), 
\end{align}
where $c'_2 := -\log c_2$.    

We prove now the first inequality of (\ref{ine_TA}). We can assume $\sup_X T_A >0$ because otherwise the desired inequality  is trivial for  any $\lambda_1 \ge 0$.  Put $\varphi_A:=T^*_A- \sup_X T^*_A$ which is an $\omega$-p.s.h. function with  $\sup_X \varphi_A =0.$ It follows that 
\begin{align}\label{ine_LpvarphiA}
\|\varphi_A\|_{L^p} \lesssim 1
\end{align}
for every $p \ge 1.$     

Let $\varphi$ be an $\omega$-p.s.h. function such that $0 \le \varphi \le 1.$ Since  $(\sup_X T_A)^{-1}\varphi_A = -1$ on $A\backslash Q_A,$  and $\capK_{BTK}(Q_A)=0,$ we obtain 
\begin{align} \label{ine_uAAA}
\int_A (\ddc \varphi+ \omega)^k \le   (\sup_X T_A)^{-1}\int_X [- \varphi_A] (\ddc \varphi+ \omega)^k \lesssim (\sup_X T_A)^{-1} \|\varphi_A\|_{L^1}
\end{align}        
for every $\varphi$ with $0 \le  \varphi \le 1$ by the Chern-Levine-Nirenberg inequality.   Combining  (\ref{ine_uAAA}) with (\ref{ine_LpvarphiA}) gives the first inequality of (\ref{ine_TA}). It remains to prove the second one. 

Recall that $-1 \le (\sup_X T_A)^{-1} \varphi_A \le 0$ and $(\sup_X T_A)^{-1} \varphi_A$ is an $(\sup_X T_A)^{-1} \omega$-p.s.h. function. Hence  $(\sup_X T_A)^{-1} \varphi_A$ is $\omega$-p.s.h. if $(\sup_X T_A)^{-1}\le 1$. Consider  the case where  $(\sup_X T_A)^{-1} \le 1.$    By definition of $\capK_{BTK},$   we get 
\begin{align}\label{ine_capBTKTA}
\capK_{BTK}(A) \ge (\sup_X T_A)^{-k} \int_A (\ddc \varphi_A+ \omega)^k=(\sup_X T_A)^{-k} \int_A (\ddc T^*_A+ \omega)^k
\end{align}
By Proposition \ref{pro_MATAuA}, we have 
$$\int_A (\ddc T^*_A+ \omega)^k= \int_X  (\ddc T^*_A+ \omega)^k.$$
Hence the second inequality of (\ref{ine_TA}) follows if $(\sup_X T_A)^{-1} \le 1.$ When $(\sup_X T_A)^{-1} \ge 1,$ then $T^*_A-1 \le 0$ on $X$ and  $\le -1$ on $A \backslash Q_A.$   We imply that
$$\capK_{BTK}(A)= \capK_{BTK}(A \backslash Q_A) \ge \int_A (\ddc T^*_A+ \omega)^k>0$$ 
which combined with the fact that $\sup_X  T_A \ge 0$ yields the second inequality of  (\ref{ine_TA}) in this case. The proof is finished.
\endproof

\begin{proof}[End of the proof of Theorem \ref{the_pluripolar}]  First observe that a countable union of pluripolar sets is again a pluripolar set. Indeed, let $(V_n)_{n\in \N}$ be a countable family of pluripolar sets on $X$. Hence we have $V_n \subset \{\varphi_n= -\infty\}$  for some  $\omega$-p.s.h function $\varphi_n$ with $\sup_X \varphi_n= 0.$  Define 
$$\varphi:= \sum_{n=1}^\infty \varphi_n/ n^2$$which is of bounded $L^1$-norm because $\|\varphi_n\|_{L^1}$ is uniformly bounded in $n.$ Hence $\varphi$ is a quasi-p.s.h. function and $V_n \subset \{\varphi= -\infty\}$ for every $n.$

Let $V$ be a locally pluripolar set.  We need to prove $V$ is pluripolar. If $V$ is compact, the desired claim is a direct application of (\ref{ine_TAorigin}). For the general case, we need some more arguments.  

By Lindel\"of's property, we can cover $V$ by at most countably many sets of form $\{\varphi_j= -\infty\}$ for some  p.s.h functions  $\varphi_j$ on some open subset $U_j$ of $X.$ Hence in order to prove the desired assertion, we only need to consider $V=\{\varphi= -\infty\}$ for some p.s.h. function $\varphi$ in an open subset $U$ of $X$ which is biholomorphic to a ball in $\C^k.$ 

  Let $U_1$ be a  relatively compact open subset of  $U.$ Suppose that $V\cap U_1$ is not pluripolar. Hence $T^*_{V \cap U_1}$ is a bounded $\omega$-p.s.h function.  Consider a decreasing sequence of smooth  p.s.h. functions  $(\varphi_n)_{n \in \N}$ defining on an open neighborhood of $\overline U_1$ converging pointwise to $\varphi.$ For every positive integer $N,$  put 
$$V_{n, N}:= \{\varphi_n \le -N\} \cap \overline U_1$$
which is a compact subset increasing in $n$. 
Hence $(T_{V_{n,N}}^*)_{n \in \N}$ is  a decreasing sequence of $\omega$-p.s.h. functions which converges pointwise to an $\omega$-p.s.h. function $T_N.$ 

Since $\{\varphi_n <-N\}$ is open, $T^*_{V_{n,N}}= T_{V_{n,N}}=0$ on $\{\varphi_n <-N\}\cap U_1.$  Thus $T_N= 0$ on $\{\varphi < -N\}\cap U_1$ which contains $V \cap U_1.$ We infer that 
\begin{align*}
0 \le T_N \le T^*_{V\cap U_1}
\end{align*}
 for every $N.$ This combined with the fact that  $(T_N)_{N \in \N}$ is increasing gives
\begin{align} \label{hoituTn}
0 \le T_\infty:= (\lim_{N \to \infty} T_N)^* \le T^*_{V \cap U_1}
\end{align}
and $T_\infty$ is an $\omega$-p.s.h. function.   Applying (\ref{ine_TAorigin}) to $A:= V_{n,N}$ we get
\begin{align}\label{ine_supTVnNMn}
\sup_X T^*_{V_{n,N}} \ge c'_2+ \lambda_2' M_{n,N}^{1/k} \capK_{BTK}(V_{n,N})^{-1/k}
\end{align}
where $M_{n,N}:= \int_X (\ddc T^*_{V_{n,N}}+ \omega)^k.$ By the convergence of Monge-Amp\`ere operators, we have 
\begin{align} \label{eq_hoituMnn}
\lim_{n \to \infty}M_{n,N}= \int_X (\ddc T_N+ \omega)^k=:M_N, \quad \lim_{N \to \infty}M_N= \int_X (\ddc T_\infty+ \omega)^k=: M_\infty 
\end{align}
Note that $M_\infty>0$ by Lemma \ref{le_boundecap}. On the other hand, we have 
 $$\capK_{BTK}(V_{n,N}) \lesssim N^{-1}$$
by the Chern-Levine-Nirenberg inequality. This together with  (\ref{eq_hoituMnn}) and   (\ref{ine_supTVnNMn}) implies 
\begin{align}\label{ine_supTVnNMn2}
\sup_X T_N \ge c'_2+ \lambda_2' M_{N}^{1/k} N^{1/k}.
\end{align}
Letting $N \to \infty$ in the last inequality and using (\ref{eq_hoituMnn}), (\ref{hoituTn}), we get 
$$ \sup_X T^*_{V \cap U_1}  \ge  \sup_X T_\infty = \infty.$$
This is a contradiction.  Hence $V\cap U_1$ is pluripolar for every relatively compact open subset  $U_1$ of $U$. It follows that  $V$ is pluripolar. This finishes the proof. 
\end{proof}

\bibliography{biblio_family_MA,bib_expose_cnrs}
\bibliographystyle{siam}

\end{document}